\documentclass[12pt,letterpaper]{article}
\usepackage{amsfonts}
\usepackage{amssymb}
\usepackage{latexsym}  % latex2e symbol package

\newtheorem{theorem}{Theorem}[section]
\newtheorem{definition}[theorem]{Definition}
\newtheorem{lemma}[theorem]{Lemma}

\newcommand{\UU}{\mathcal{U}}

\newcommand{\HS}{\mathit{HS}}
\newcommand{\MA}{\mathit{MA}}
\newcommand{\CH}{\mathit{CH}}
\newcommand{\ZFC}{\mathit{ZFC}}
\newcommand{\PFA}{\mathit{PFA}}
\newcommand{\Iseq}{\mathit{Iseq}}

\newcommand{\cl}{\mathrm{cl}}

\newcommand\res{\mathord {\upharpoonright}}  % less space around it

\newenvironment{itemizz}{\begin{itemize}\setlength{\itemsep}{-1mm}} %
{\end{itemize}}

\newcommand{\eop}{$\quad\spadesuit $}
\newenvironment{proof}{{\bf Proof.}}{\eop\medskip}

\begin{document}

\title{
A Compact Homogeneous S-space}

\author{Ramiro de la Vega\footnote{University of Wisconsin,  Madison, WI  53706, U.S.A.,
\ \ delavega@math.wisc.edu}
\  and
Kenneth Kunen\footnote{University of Wisconsin,  Madison, WI  53706, U.S.A.,
\ \ kunen@math.wisc.edu; partially supported by NSF Grant DMS-0097881.}
}

\maketitle

\begin{abstract}
Under the continuum hypothesis, there is a compact homogeneous strong
S-space.
\end{abstract}

\section{Introduction}

A space $X$ is \textit{hereditarily separable} ($\HS$) iff every subspace is
separable. An $S$\textit{-space} is a regular Hausdorff $\HS$ space with a
non-Lindel\"of subspace. A space $X$ is \textit{homogeneous} iff for every
$x,y\in X$ there is a homeomorphism $f$ of $X$ onto $X$ with
$f( x) =y$.  Under $\CH$,
several examples of $S$-spaces have been constructed, including
topological groups (see \cite{hj}) and compact $S$-spaces (see \cite{jkr}).
It is asked in \cite{arc,kph} whether there are compact homogeneous $S$-spaces.
As we shall show in Theorem \ref{thm-homsp}, there are under $\CH$.
This cannot be done in $\ZFC$, since
there are no compact $S$-spaces
under $\MA+\neg \CH$ (see \cite{sze}); 
there are no $S$-spaces at all under $\PFA$ (see \cite{tod}).

In Section \ref{sspace}, we use a slightly modified version of the
construction in \cite{jkr,neg} to refine the topology of any given second
countable space, and turn it into a first countable
\textit{strong $S$-space} (i.e.,
each of its finite powers is an $S$-space). In Section \ref{comp}, we show
that if the original space is compact, then there is a natural
compactification of the new space which is also a first
countable strong S-space.  If in addition the
original space is zero-dimensional, then
the $\omega^\mathrm{th}$ power of this compactification
will be homogeneous by Motorov \cite{mot}, proving Theorem \ref{thm-homsp}.

\section{A Strong S-Space
\label{sspace}}

If $\tau $ is a topology on $X$,
we write $\tau ^I$ for the corresponding product
topology on $X^I$; likewise if $\tau ^{\prime }\subseteq \tau $ is a base we
write $( \tau ^{\prime }) ^I$ for the natural corresponding base
for $\tau ^I$.
If $E \subseteq X$, then $\cl(E, \tau)$ denotes the closure of $E$
with respect to the topology $\tau$.
This notation will be used when we are discussing
two different topologies on the same set $X$.

The following two lemmas are well-known; the second
is Lemma 7.2 in \cite{neg}:

\begin{lemma}
\label{lemma-hs2}
If $X$ is $\HS$ and $Y$ is second countable, then $X\times Y$ is $\HS$
\end{lemma}

\begin{lemma}
\label{lemma-hsomega}
$X^\omega$ is $\HS$ iff $X^n$ is $\HS$ for all $n<\omega$.
\end{lemma}

The next lemma, an easy exercise, is used in the proof of
Theorem \ref{thm-strong-ssp}:

\begin{lemma}
\label{lemma-cl-prod}
If $(x,y) \in X \times Y$ and $S \subseteq X \times Y$,
then $(x,y) \in \cl(S)$ iff
$y \in \cl(\pi(S \cap (U \times Y)))$ for all neighborhoods $U$ of $x$,
where $\pi : X \times Y \to Y$ is projection.
\end{lemma}

The following is proved (essentially) in \cite{neg}, but 
our proof below may be a bit simpler:

\begin{theorem}
\label{thm-strong-ssp}
Assume $\CH$.  Let  $\rho$ be a second countable $T_3$
topology on $X$, where $|X| = \aleph_1$.  Then there is a finer topology $\tau$
on $X$ such that $(\omega_1, \tau)$ is a first countable
locally compact strong S-space.
\end{theorem}
\begin{proof}
WLOG, $X = \omega_1$.
For $\eta <\omega _1$ we write $\rho _\eta $ for the topology of $\eta $ as a
subspace of $\left( \omega_1,\rho \right)$.
Applying $\CH$, list 
$\bigcup_{0<n<\omega }[ (\omega_1)^n]^{\leq \omega }$ as
$\{ S_\mu :\mu \in \omega _1\}$, so that each 
$S_\mu \subseteq \mu ^{n( \mu )}$ for some $n(\mu)$ with 
$0 < n(\mu) < \omega$.

For $\eta \leq \omega _1$ we construct $\tau _\eta $ a topology on $\eta $
by induction on $\eta $ so as to make the following hold for all $\xi <\eta
\leq \omega _1$:

\begin{itemizz}
\item[1.]
$\tau _\xi =\tau _\eta \cap {\cal P}\left( \xi \right) $.
\item[2.]
$\tau _\eta $ is first countable, locally compact, and $T_3$.
\item[3.]
$\tau _\eta \supseteq \rho _\eta $.
\end{itemizz}

Note that (1) implies in particular that $\xi \in \tau_\eta$;
that is, $\xi$ is open.  Thus, if $\tau = \tau_{\omega_1}$, then
$(\omega_1, \tau)$ is not Lindel\"of.  Also by (1),
$\tau_\eta$ for limit $\eta$ is determined from the $\tau_\xi$
for $\xi < \eta$.
So, we need only specify what happens at successor ordinals.

For $n \ge 1$ and $\xi < \omega_1$,
let $\Iseq(n, \xi)$ be the set of all $f \in (\omega_1)^n$
which satisfy $f(0) < f(1) < \cdots < f(n-1) = \xi$.
The following condition states
our requirement on $\tau_{\xi+1}$:

\begin{itemizz}
\item[4.] For each $\mu < \xi$ and each $f \in \Iseq(n, \xi)$,
where $n = n_\mu$:
$$
f \in \cl(S_\mu,\, (\tau_{\xi+1})^{n-1}\times \rho) 
\ \Longrightarrow \ f \in \cl(S_\mu,\, (\tau_{\xi+1})^{n}) \ \ .
$$
\end{itemizz}
If $n = n_\mu = 1$, then $(\tau_{\xi+1})^{n-1}\times \rho$ just
denotes $\rho$.  That is, (4) requires
$$
\xi \in \cl(E,\rho) \ \Longrightarrow \  \xi \in \cl(E,\tau_{\xi+1})
\eqno{(*)}
$$
for all $E$ in the countable family
$\{S_\mu : \mu < \xi \ \&\ n(\mu) = 1\}$.
It is standard (see \cite{jkr}) that
one may define $\tau_{\xi+1}$ so that this holds.
Now, consider (4) in the case $n = n_\mu \ge 2$.
By (2), $\tau _\xi $ is second countable, so let $\tau _\xi'$ be
a countable base for $\tau_\xi$.  Applying Lemma \ref{lemma-cl-prod},
(4) will hold if whenever 
$U = U_0 \times \cdots \times U_{n-2} \in (\tau_\xi')^{n-1}$
is a neighborhood of $f \res (n-1)$, 
$$
\xi \in \cl(\pi(S_\mu \cap (U \times (\xi+1))),\rho) \ \Longrightarrow \  
\xi \in \cl(\pi(S_\mu \cap (U \times (\xi+1))),\tau_{\xi+1}) \ \ ,
$$
where $\pi : \xi^{n-1} \times (\xi + 1) \to (\xi + 1)$ is projection.
But this is just a requirement of the form $(*)$ for countably many
more sets $E$, so again there is no problem meeting it.

Now, we need to show that $\tau ^n$ is $\HS$ for each $0<n<\omega $.
We proceed by induction, so assume that  $\tau ^m$ is $\HS$
for all $m < n$.  Fix $A \subseteq (\omega_1)^n$; we need to show that
$A$ is $\tau^n$-separable.
Applying the induction hypothesis, we may assume that each $f\in A$
has all coordinates distinct.  Also, since permutation of coordinates
induces a homeomorphism of $(\omega_1)^n$, we may assume that
each $f\in A$ is strictly increasing; that is,
$f \in \Iseq(n, \xi)$, where $\xi = f(n-1)$.
By the induction hypothesis and Lemma \ref{lemma-hs2}, $A$
is separable in $(\tau_{\xi+1})^{n-1}\times \rho$.
We can then fix $\mu$ such that $n(\mu) = n$, $S_\mu \subseteq A$,
and $S_\mu$ is $(\tau_{\xi+1})^{n-1}\times \rho$-dense in $A$.
Now, say $f\in A$ with $\xi =  f(n-1) > \mu$.  Applying (4),
we have $f \in \cl(S_\mu,\, \tau^{n})$.
Thus, $A \setminus \cl(S_\mu,\, \tau^{n})$ is countable,
so $A$ is $\tau^n$-separable.
\end{proof}

\section{Compactification\label{comp}}

\begin{definition}
\label{def-dotcup}
If $\varphi$ is a continuous map from the $T_2$ space  $Y$ into $X$, then
$Y \dot\cup_\varphi X$ denotes the disjoint union of $X$ and $Y$,
given the topology which has as a base:
\begin{itemizz}
\item[a.] All open subsets of $Y$, together with
\item[b.] All $[U,K] := U\cup (\varphi ^{-1}U \setminus K)$,
where $U$ is open in $X$ and $K$ is compact in $Y$.
\end{itemizz}
\end{definition}

Our main interest here is in the case where $X$ is compact and $Y$
is locally compact.  Then, if  $|X| = 1$, we have the 1-point 
compactification of $Y$, and if $Y$ is discrete and $\varphi$ is a bijection
we have the Aleksandrov duplicate of $X$.

\begin{lemma}
\label{lemma-properties-dotcup}
Let $Z = Y \dot\cup_\varphi X$, with $X,Y$ Hausdorff:
\begin{itemizz}
\item[1.] $X$ is closed in $Z$, $Y$ is open in $Z$, and both
$X,Y$ inherit their original topology as subspaces of $Z$.
\item[2.] If $Y$ is locally compact, then $Z$ is Hausdorff.
\item[3.] If $X$ is compact, then $Z$ is compact.
\item[4.] If $X,Y$ are first countable, $X$ is compact,
$Y$ is locally compact, and each $\varphi^{-1}(x)$ is compact,
then $Z$ is first countable.
\item[5.] If $X,Y$ are zero dimensional, $X$ is compact,
and $Y$ is locally compact, then $Z$ is zero dimensional.
\item[6.] If $X$ is second countable and $Y^\omega$ is $\HS$,
then $Z^\omega$ is $\HS$.
\end{itemizz}
\end{lemma}
\begin{proof}
For (3):
If $\UU $ is a basic open cover of $Z$, then there are
$n\in \omega $ and $[ U_i,K_i] \in \UU $ for $i < n$ such that
$\bigcup_{i < n}U_i=X$.  Thus,
$\bigcup_{i < n}[U_i,K_i]$ contains all points of $Z$ except for
(possibly) the points in the compact set $\bigcup_{i < n}K_i \subseteq Y$.

For (4): $Z$ is compact Hausdorff and of countable pseudocharacter.

For (5): $Z$ is compact Hausdorff and totally disconnected.

For (6): By Lemma \ref{lemma-hsomega}, it is sufficient to
prove that each $Z^n$ is $\HS$.  But $Z^n$ is a finite union
of subspaces of the form $X^j \times Y^k$, which are $\HS$
by Lemma \ref{lemma-hs2}.
\end{proof}

\section{Homogeneity}
\label{sec-conc}

The following was proved by Dow and Pearl \cite{dp}:

\begin{theorem}
\label{thm-dp}
If $Z$ is first countable and zero dimensional, then $Z^\omega$
is homogeneous.
\end{theorem}

Actually, we only need here the special case of this result
where $Z$ is compact and has a dense set of isolated points;
this was announced (without proof) earlier by Motorov \cite{mot}.

Note that by \v Sapirovski\u\i \ \cite{sap},
any compact $\HS$ space must have countable $\pi $-weight
(see also \cite{hod}, Theorem 7.14),
so if it is also homogeneous, it
must have size at most $2^{\aleph _0}$ by van Douwen \cite{van}. Under $\CH$
this implies, by the \v Cech -- Posp\'\i\v sil Theorem,
that the space must be first countable.

\begin{theorem}
\label{thm-homsp}
$\left( \CH\right) $ There is a (necessarily first countable)
zero-dimensional compact homogeneous strong $S$-space.
\end{theorem}
\begin{proof}
Let $X$ be the Cantor set $2^\omega $ with its usual topology,
let $Y$ be $2^\omega $ with the topology constructed in Theorem
\ref{thm-strong-ssp}, let $\varphi $ be the identity,
and let $Z = Y \dot\cup_\varphi X$.  
By Lemma \ref{lemma-properties-dotcup},
$Z$, and hence also $Z^\omega$, are
zero-dimensional first countable compact strong S-spaces;
$Z^\omega$ is homogeneous by Theorem \ref{thm-dp}.
\end{proof}

No compact topological group can be an $S$-space or an $L$-space.
However under $\CH$ there are,
by \cite{kun}, compact $L$-spaces which are right topological groups (i.e.
they admit a group operation such that multiplication on the right by a
fixed element defines a continuous map). We do not know whether there can be
compact $S$-spaces which are right topological groups.

\end{document}